\documentclass
{amsart}
\usepackage{amssymb}   
\title[Wahl's conjecture in positive characteristic]%
{Wahl's conjecture holds in odd characteristics\\%
for symplectic and orthogonal Grassmannians}
\author{V.~Lakshmibai}
	\address{Northeastern University, Boston, USA}
	\email{lakshmibai@neu.edu}
\author{K.~N.~Raghavan}
	\address{Institute of Mathematical Sciences, C.~I.~T.~Campus\newline
		\hbox{\hspace{\parindent}}%
		Chennai 600\,113, INDIA}
	\email{knr@imsc.res.in}
\author{P.~Sankaran}
	\address{Institute of Mathematical Sciences, C.~I.~T.~Campus\newline
		\hbox{\hspace{\parindent}}%
		Chennai 600\,113, INDIA}
	\email{sankaran@imsc.res.in}
\thanks{\noindent
V.~Lakshmibai was partially supported by NSF grant DMS-0652386 and
Northeastern University RSDF 07--08.
K.~N.~Raghavan and P.~Sankaran were partially supported by 
DAE grant No.~11--R\&D--IMS--5.01--0500.
This work was partly done in April~2006 during a visit by all three authors
to the {\sffamily Abdus Salam International Centre for Theoretical Physics}
whose hospitality and patronage are gratefully acknowledged.}

\usepackage{ifthen}     

\newcommand\finalized{yes}

\newcommand\ignore[1]{}

\usepackage{color}
\providecommand\wantcolor{yes}   %
\ifthenelse{\equal{\wantcolor}{no}}{\renewcommand\color[1]{}}{}
\definecolor{backgroundyellow}{cmyk}{.2,.1,.8,.2}
\definecolor{backgroundblue}{rgb}{0,0,1}
\definecolor{backgroundred}{rgb}{1,0,0}
\definecolor{backgroundmagenta}{cmyk}{0,1,0,0}
\definecolor{goodforinversevideo}{rgb}{.65,.8,1}

\newcommand\mysection{\section}

\newcommand\mysubsection{\subsection}
\newcommand\mysubsubsection[1]{%
		\subsubsection{\sffamily\upshape\mdseries #1}}
\newcommand\mysss{\mysubsubsection}

\newtheorem{annotation}{Annotation}
\newtheorem{theorem}[annotation]{
		Theorem}

\newtheorem{definition}[annotation]{
		Definition}
\newcommand\bdefn{\begin{definition}\begin{upshape}}
\newcommand\edefn{\end{upshape}\end{definition}}

\newtheorem{proposition}[annotation]{
		Proposition}
\newtheorem{example}[annotation]{
		Example}
\newcommand\bexample{\begin{example}\begin{rm}}
\newcommand\eexample{\end{rm}\hfill$\Box$\end{example}}
\newtheorem{conjecture}[annotation]{
		Conjecture}
\newcommand\bconjecture{\begin{conjecture}\begin{rm}}
\newcommand\econjecture{\end{rm}\end{conjecture}}
\newtheorem{notation}[annotation]{
		Notation}
\newcommand\bnotation{\begin{notation}\begin{rm}}
\newcommand\enotation{\end{rm}\end{notation}}

\newtheorem{remark}[annotation]{
		Remark}
\newcommand\bremark{\begin{remark}\begin{sffamily}\begin{upshape}}
\newcommand\eremark{\end{upshape}\end{sffamily}\end{remark}}
\newenvironment{myproof}{%
\par\noindent{\scshape Proof:}\begin{rm}}{\hfill$\Box$\end{rm}\newline}

\providecommand\finalized{no}
\ifthenelse{\equal{\finalized}{yes}}{}{\marginparwidth=2cm}
\ifthenelse{\equal{\finalized}{yes}}%
{\newcommand\mylabel[1]{\label{#1}}}%
{\newcommand\mylabel[1]{\label{#1}\marginpar{[{\ttfamily\upshape\tiny #1}]}}}
\ifthenelse{\equal{\finalized}{yes}}%
{\newcommand\checked[1]{}}%
{\newcommand\checked[1]{\marginpar{[{\ttfamily\upshape\tiny CHECKED: #1}]}}}
\ifthenelse{\equal{\finalized}{yes}}%
{\newcommand\spellchecked[1]{}}%
{\newcommand\spellchecked[1]{\marginpar{[{\ttfamily\upshape\tiny SPELLCHECKED: #1}]}}}

\providecommand\version{public}   
\ifthenelse{\equal{\version}{public}}%
{\newcommand\mcomment[1]{}}%
{\newcommand\mcomment[1]{\marginpar{{\sffamily\upshape\tiny #1}}}}
\ifthenelse{\equal{\version}{public}}%
{\newcommand\fcomment[1]{}}%
{\newcommand\fcomment[1]{\footnote{#1}}}

\newcommand\tensor\otimes
\newcommand\idealsheaf{\mathcal{I}}
\newcommand\idealy{\idealsheaf_Y}
\newcommand\sheafo{\mathcal{O}}
\newcommand\ox{{\sheafo_X}}
\newcommand\hash{\#}
\newcommand\Hom{\textup{Hom}}
\newcommand\hnot{H^0}
\newcommand\sheafhom{{\mathcal{H}\textup{om}}}
\newcommand\lbl{\mathcal{L}}

\newcommand\ztilde{\tilde{Z}}

\newcommand\oztilde{\sheafo_{\ztilde}}
\newcommand\ordysigma{\textup{ord}_Y\sigma}
\newcommand\sigmatilde{\tilde{\sigma}}
\newcommand\oy{{\sheafo_Y}}
\newcommand\fstar{F_*}

\newcommand\field{\pmb{k}}

\newcommand\fund\varpi
\newcommand\character\varepsilon
\newcommand\sectionp{\mathfrak{p}}
\newcommand\sectionq{\mathfrak{q}}
\newcommand\sectionr{\mathfrak{r}}

\newcommand\lbundle{\mathcal{L}}
\newcommand\lieg{\mathfrak{g}}
\newcommand\lieb{\mathfrak{b}}

\newcommand\sep{\,|\,}

\newcommand\steinberg{\textup{St}}

\newcommand\weylinvolution{\iota}

\newcommand\End{\textrm{End}}
\newcommand\basis{e}

\usepackage{ragged2e}    
\begin{document}
\begin{abstract}%
It is shown that the proof by Mehta and Parameswaran of Wahl's conjecture
for Grassmannians in positive odd characteristics
also works for symplectic and orthogonal Grassmannians.
\end{abstract}

\maketitle
Let $X$ be a non-singular projective variety over an algebraically closed
field~$\field$.     For ample line bundles $L$ and $M$ over $X$,
consider the natural restriction map (called the {\em Gaussian\/})
\begin{equation}\label{e.gaussian}
 H^0(X\times X, \idealsheaf_\Delta\tensor p_1^*L\tensor p_2^*M)\to
H^0(X, \Omega^1_{X/\field}\tensor L\tensor M)
\end{equation} 
where $\idealsheaf_\Delta$ denotes the
ideal sheaf of the diagonal $\Delta$ in $X\times X$, $p_1$ and $p_2$
the two projections of $X\times X$ on $X$, and $\Omega^1_{X/\field}$
the sheaf of differential $1$-forms of~$X$ over~$\field$.    
Wahl conjectured in~\cite{wahl} 
that this map is surjective when $X$ is a homogeneous space for the action
of a semisimple group $G$,  that is,  when $X=G/P$ for 
$G$ a semisimple and simply connected linear algebraic group over~$\field$
and  $P$ a parabolic subgroup of $G$.    The original conjecture
was perhaps meant only over the field of complex numbers,
and in fact it has been proved in that case by Kumar~\cite{kumar} 
using representation theoretic techniques,
but following~\cite{lmp,mp} we use the term `Wahl's conjecture' to refer
to the surjectivity of the Gaussian without any restriction on the
characteristic.    The truth
of the conjecture in infinitely many positive characteristics
would imply its truth in characteristic zero, for the Gaussian 
is defined over the integers.

Assume from now on that the base field $\field$ has positive characteristic.
Lakshmibai, Mehta, and Parameswaran~\cite{lmp} show that,
in odd characteristic, 
Wahl's conjecture holds if there is a Frobenius splitting of $X\times X$ that 
compatibly splits the diagonal~$\Delta$ and has maximal multiplicity
along $\Delta$ (the definitions are recalled below). 
Moreover they conjecture that such a Frobenius splitting exists
(in any characteristic).    
Mehta and Parameswaran~\cite{mp} prove that this latter conjecture
holds for Grassmannians. In the present paper
it is shown that their proof also works
for symplectic and orthogonal Grassmannians---see Theorem~\ref{t.main}
and the conclusion in~\S\ref{s.conclusion} 
below.

This paper is organized as follows:
notation is fixed in~\S\ref{s.notation},  some basic definitions and
results about Frobenius splittings are recalled in~\S\ref{s.splittings},
the results of~\cite{lmp} about splittings for blow-ups 
are recalled in~\S\ref{s.blowups},   
some important results about splittings
for complete homogeneous spaces are recalled in~\S\ref{s.canonical},
the main result (Theorem~\ref{t.main}) is 
proved in~\S\ref{s.theorem},  and Wahl's conjecture for ordinary,
symplectic, and orthogonal Grassmannians is deduced from the main theorem
in~\S\ref{s.conclusion} following the argument in~\cite{mp}.
\mysection{Notation}\mylabel{s.notation}
The following notation remains fixed throughout:
\begin{itemize}
\item $\field$ an algebraically closed field of positive characteristic~$p$;
\item $V^{[m]}$  where $V$ is a $\field$-vector space and $m$ an integer
denotes the $\field$-vector space obtained by pulling back~$V$ via
the automorphism $\lambda\mapsto\lambda^{p^m}$ of~$\field$.
\item $G$ a semisimple, simply connected group linear algebraic over $\field$;
\item $B$ a Borel subgroup of $G$;
\item $T$ a maximal torus of $G$ such that $T\subseteq B$;
\item $P$ a {\em standard\/} parabolic subgroup of $G$ (standard means
  $B\subseteq P$);
\item the {\em roots\/} are taken with respect to~$T$;
\item the {\em positive roots\/} are taken to be the characters for the
adjoint action of~$T$ on the Lie algebra of the unipotent part~$U$ of~$B$;
\item $\rho$ denotes half the sum of all positive roots (equivalently,
the sum of all fundamental weights);
\item $w_0$ denotes the longest element of the Weyl group;
\item $\weylinvolution$ denotes the Weyl involution $\lambda\mapsto
-w_0\lambda$ on characters of~$T$;
\item for a character~$\lambda$ of~$T$,  
\begin{itemize}
\item
$\lambda$ denotes also
its extension to~$B$ via the isomorphism $B/U\cong T$ (induced by the inclusion
of~$T$ in~$B$);
\item   
$\field_{\lambda}$ denotes the one dimensional $B$-module defined by the character $\lambda$;
\item $\lbundle(\lambda)$ denotes the line bundle $G\times_B k_{-\lambda}$
on~$G/B$;
\item $H^0(G/B,\lbundle(\lambda))$ denotes the $G$-module of global sections
of $\lbundle(\lambda)$.
\end{itemize}
\end{itemize}
Observe that
$H^0(G/B,\lbundle(\lambda))$ can be identified as
a $G$-module with the space of regular functions~$f$ on~$G$ that
transform thus:
\begin{equation}
  \label{eq:transform1}
  f(gb)=\lambda(b)f(g) \quad \forall~g\in G \quad \forall~b\in B
\end{equation}
the action of $G$ on functions being given by $(gf)(x):=f(g^{-1}x)$;
the highest and lowest weights of~$H^0(G/B,\lbundle(\lambda))$ are
respectively $-w_0\lambda$ and $-\lambda$.
Observe also that the anti-canonical 
bundle~$K^{-1}$ of $G/B$ is~$\lbundle(2\rho)$,
for $G\times^B{\lieg/\lieb}\to G/B$ is the tangent bundle of $G/B$,
where $\lieg$ and $\lieb$ denote the Lie algebras of $G$ and $B$ respectively.

\mysection{Frobenius Splittings}\mylabel{s.splittings}
Let $X$ be a scheme over~$\field$, separated and of finite type.
Denote by $F$ the {\em absolute Frobenius\/} map on $X$: this is
the identity map on the underlying topological space~$X$ and is
the $p$-th power map on the structure sheaf $\ox$.    
We say that $X$ is {\em Frobenius split\/} if the $p$-th power map
$F^\hash:F_*\ox\leftarrow\ox$ splits as a map of $\ox$-modules
(see~\cite[\S1,~Definition~2]{mr}, \cite[Definition~1.1.3]{bk}).
A splitting $\sigma: F_*\ox\to\ox$ {\em compatibly splits\/} 
a closed subscheme~$Y$
of $X$ if $\sigma(F_*\idealy)\subseteq \idealy$ where $\idealy$ is the
ideal sheaf of~$Y$ 
(see \cite[\S1,~Definition~3]{mr}, \cite[Definition~1.1.3]{bk}).

Now suppose that $X$ is a non-singular projective variety and
denote by~$K$ its canonical bundle.
Using Serre duality and the observation that $F^*\lbl\cong\lbl^p$ for
an invertible sheaf~$\lbl$ on $X$ (applied to~$K$),  we get a canonical
isomorphism of $\hnot(X,\sheafhom(F_*\ox,\ox))=\Hom_\ox(F_*\ox,\ox)$
with $\hnot(X,K^{1-p})^{[1]}$
(see~\cite[Page~32]{mr} or \cite[Lemma~1.2.6 and \S1.3]{bk}; 
without the superscript `$[1]$',  
the isomorphism would only be $\field$-semilinear).
To say that {\em $\sigma$ splits $X$\/} for 
$\sigma\in\hnot(X,K^{1-p})\cong\Hom_\ox(\fstar\ox,\ox)^{[-1]}$ means
that the underlying~$\ox$-module homomorphism $F_*\ox\to\ox$ is a
splitting of~$F^\hash$.     
Set \[ \boxed{\End_F(X):=\Hom_\ox(\fstar\ox,\ox)^{[-1]}}   \]


\mysection{Splittings and Blow-ups}\mylabel{s.blowups}
Let $Z$ be a non-singular projective variety and $\sigma$ a section
of $K^{1-p}$ (where $K$ is the canonical bundle) that splits~$Z$.
Let~$Y$ be a closed non-singular subvariety 
of~$Z$ of codimension~$c$.   Let $\ordysigma$
denote the order of vanishing of $\sigma$ along $Y$. Let $\ztilde$ denote
the blow up of $Z$ along~$Y$ and~$E$ the exceptional divisor (the fiber
over~$Y$) in~$\ztilde$.    

A splitting $\tilde{\tau}$ of
$\ztilde$ induces a splitting~$\tau$ on~$Z$,
by virtue of\/ $\pi_*\mathcal{O}_{\ztilde}\leftarrow\mathcal{O}_Z$ being an isomorphism
where $\pi:\ztilde\to Z$ is the natural map (see the result to this effect
quoted in \S\ref{s.canonical} below).
We say that~$\sigma$ {\em lifts\/} to a
splitting of $\ztilde$ if it is induced thus from a splitting~$\tilde{\sigma}$
of~$\ztilde$.   The lift of $\sigma$ to $\ztilde$ is unique if it exists,
since $\ztilde\to Z$ is birational and two global sections of the
locally free sheaf $\sheafhom_{\oztilde}(\fstar\oztilde,\oztilde)$ 
that agree on an open set must be equal.

\begin{proposition}
With notation as above, we have 
\begin{enumerate}
\item $\ordysigma\leq c(p-1)$.
\item
If $\ordysigma=c(p-1)$ then
$Y$ is compatibly split.
\item $\ordysigma\geq(c-1)(p-1)$ if and only if $\sigma$ lifts to a
splitting~$\sigmatilde$ of~$\ztilde$; moreover, 
$\ordysigma\geq c(p-1)$ if and only if the splitting
$\sigmatilde$ is compatible with $E$.
\end{enumerate}
\end{proposition}
\begin{myproof}
These statements appear as Exercise~1.3.E.12 in~\cite{bk}.   In any
case,~(1) and~(2) are elementary to see from the local description 
as in~\cite[Proposition~5]{mr}
of the functorial isomorphism between~$\End_F(X)$
and $H^0(X, K^{1-p})$.
Statement~(3) is Proposition~2.1 of~\cite{lmp}.
\end{myproof}

\noindent
The proposition above justifies the definition below:
\bdefn\textup{(\cite[Remark~2.3]{lmp})}
We say that $Y$ is {\em compatibly split by $\sigma$ with maximal
multiplicity\/} if $\ordysigma=c(p-1)$.
\edefn

Now let $Z=G/P\times G/P$, and $Y$ the diagonal copy of $G/P$ in~$Z$.
We have:
\begin{theorem}\label{t.lmp}
  \textup{(
\cite[pages~106--7]{lmp})}
Assume that the characteristic~$p$ is odd.
If $E$ is compatibly split in $\ztilde$, or, equivalently,  if there is
a splitting of $Z$ compatibly splitting~$Y$ with maximal multiplicity,
then the Gaussion map~(\ref{e.gaussian}) is surjective for~$X=G/P$.
\end{theorem}

\bconjecture
  \textup{(
\cite[page~106]{lmp})}  For any $G/P$,
there exists a splitting of $Z$ that compatibly splits
the diagonal copy of~$G/P$ with maximal multiplicity.
\econjecture

\mysection{Canonical splitting}\mylabel{s.canonical}
For a $B$-scheme~$X$ there is the notion of a 
{\em $B$-canonical element\/} in
$\End_F(G/P)$ 
(\cite[Definition~4.3.5]{vdk}, \cite[Definitions~4.1.1,~4.1.4]{bk}).
We can take the following characterization to
be the definition:
\begin{proposition}\mylabel{p:canonical}
{\begin{rm}(\cite[Lemma~4.1.6]{bk})\end{rm}}
For a $B$-scheme $X$,   an element~$\phi$ belonging to
 $\End_F(X)$ is 
$B$-canonical if and only if there exists a $\field$-linear $B$-module map
\[
\Theta_\phi:\steinberg\tensor \field_{(p-1)\rho}\to\End_F(X)
\quad\textup{with $\Theta_\phi(f_-\tensor f_+)=\phi$} \]
where $\field_{(p-1)\rho}$ is the one-dimensional $B$-module on which $B$ acts
by the character $(p-1)\rho$, $f_-$ is a non-zero lowest weight vector
of the Steinberg module $\steinberg:=H^0(G/B,\lbundle((p-1)\rho)$,
and $0\neq f_+\in \field_{(p-1)\rho}$.
\end{proposition}
%
\begin{theorem}\mylabel{t:canonical} \textup{(\cite{mat}; \cite[Theorem~4.1.15]{bk})}  There is a unique
(up to non-zero scalar multiples) non-zero $B$-canonical element in
$\End_F(G/P)$. 
Moreover,  this element is a 
splitting of~$G/P$
and compatibly splits all Schubert and opposite Schubert subvarieties.
\end{theorem}
\noindent
Since, on the one hand, splittings of $X$ are mapped to splittings
of $Y$ under $f_*$ for 
any morphism $f:X\to Y$ of schemes such that $f_*\ox\leftarrow\oy$
is an isomorphism (\cite[Proposition~4]{mr}, \cite[Proposition~1.1.8]{bk}),
and, on the other hand, as is readily seen 
(see also~\cite[Exercise~4.1.E.3]{bk}),
$B$-canonical elements of $\End_F(Y)$ 
are mapped to $B$-canonical elements of~$\End_F(X)$ 
by $f_*$ for such~$f$ that are $B$-morphisms of $B$-schemes,   
it follows that the $B$-canonical
splitting of $G/B$ is mapped to the $B$-canonical splitting of $G/P$ under the
natural map $G/B\to G/P$.
\begin{theorem}\mylabel{t.canonical.2}
\textup{(\cite[Theorem~4.1.17, Remark~4.1.18]{bk}, \cite{mat00})}
For a $B$-scheme $X$ there is a natural injective association
$\sigma\mapsto\tilde{\sigma}$ of 
$B$-canonical elements in $\End_FX$ to $B$-canonical elements of 
$\End_F(G\times^B X)$.   Splittings are mapped to splittings under
this association.   Moreover, for a $B$-canonical splitting $\sigma$ of~$X$,
the splitting~$\tilde{\sigma}$ is the unique 
$B$-canonical splitting of $G\times^B X$ that
compatibly splits $X\cong e\times X\subseteq G\times^B X$ (the fiber over
the identity coset of~$G/B$) and restricts on~$e\times X$ to~$\sigma$.
\end{theorem}
\noindent
Consider the isomorphism $G\times^B G/B\cong G/B\times G/B$ defined 
by the association 
$(g_1,g_2B)\mapsto(g_1B,g_1g_2B)$.    It follows from the above
theorems that there exists a unique $B$-canonical splitting of $G/B\times G/B$
that compatibly splits~$e\times G/B$ and restricts to the canonical splitting
of $G/B\cong e\times G/B$---here $G/B\times G/B$ is a $B$-variety
by the diagonal action.    This we call 
{\em the canonical splitting\/} of~$G/B\times G/B$.
The splitting of $G/P\times G/P$ obtained as the push forward of this
under the natural map $G/B\times G/B\to G/P\times G/P$ is called
{\em the canonical splitting\/} of $G/P\times G/P$.

\mysection{The theorem}\mylabel{s.theorem}  
\begin{theorem}\mylabel{t.main}
The $B$-canonical splitting of $G/B$ has maximal multiplicity along~$P/B$
in the following cases (we are fixing a maximal torus~$T$ of~$G$ contained
in~$B$, and the ordering of the simple
roots (and so also the fundamental weights) is as in 
Bourbaki~\cite{bour}):
\begin{enumerate}
\item $G=SL_n$ and $P$ is any maximal parabolic (the set of such
$G/P$ are precisely all Grassmannians).  
\textup{(This is already in~\cite{mp}, but we prove it again below.)}
\item $G=Sp_{2n}$ and $P=P_n$ (the set of such $G/P$ are precisely
all symplectic Grassmannians).
\item
the characteristic is $\geq3$, $G=SO_{2n}$ and $P=P_n$ 
(the set of such $G/P$ are precisely all orthogonal Grassmannians: as
is well known $SO_{2n+1}/P_n\cong SO_{2n+2}/P_{n+1}$;  but it is not
true that the $B$-canonical splitting of~$SO_{2n+1}/B$ has maximal
multiplicity along $P_n/B$).%
\footnote{Our convention that $G$ be simply connected is
violated here.  The violation is however not serious and 
we trust the reader can make the appropriate adjustments.}

\end{enumerate}
\end{theorem}
\begin{myproof} 
%
Let $\fund_1$, \ldots, $\fund_\ell$ be the fundamental
weights ordered as in Bourbaki, and~$\rho$ their sum
$\fund_1+\cdots+\fund_\ell$.    
Global sections of $K^{-1}$ can be identified with regular functions on~$G$ that
transform thus:
\begin{equation}\label{eq:transform}
{f(gb)} = \rho(b)^{2} f(g)\quad \forall\ g\in G,\quad \forall\ b\in B
\end{equation}
We will explicitly write, as such a function,  the section of $K^{-1}$
whose $(p-1)^{\textup{st}}$ power gives the $B$-canonical splitting of~$G/B$.

\mysubsection{The case $G=SL_n$ and $P$ any maximal parabolic}%
\mylabel{ss.sln}
Although this case is dealt with
already by Mehta and Parameswaran~\cite{mp},  rewriting their proof 
as below is helpful.   
We take the Borel subgroup~$B$ to be
the subgroup consisting of upper triangular
matrices, and the maximal torus~$T$ to be the subgroup consisting
of diagonal matrices.
We denote by $\character_k$ the character of~$T$ that maps elements
of $T$ to their $(k,k)$-entries.  The simple roots in order are
$\character_1-\character_2$, $\character_2-\character_3$, \ldots,
$\character_{n-1}-\character_{n}$, and the corresponding
fundamental weights are 
\[ \fund_1=\character_1,\quad \fund_2=\character_1+\character_2,\quad
\ldots,\quad \fund_{n-1}=\character_1+\cdots+\character_{n-1}.  \]

For $1\leq k\leq n-1$ and $1\leq a_1<\ldots<a_k\leq n$,  let 
$\det(a_1,\ldots,a_k\sep1,\ldots,k)$ denote the function on the space of 
$n\times n$ matrices (and so also on the group $SL_n$ by restriction) 
obtained by taking the determinant of the submatrix
formed by the first~$k$ columns and the rows numbered $a_1$, \ldots, $a_k$.
\begin{proposition}\mylabel{p.determinants}
\begin{enumerate}
\item \mylabel{i.one}
$\det(a_1,\ldots,a_k\sep1,\ldots,k)$ is a global section of the line bundle
$\lbundle(\fund_k)$; it is
a weight vector of weight~$-(\character_{a_1}+\cdots+\character_{a_k})$.
\item\mylabel{i.two}
$\det(1,\ldots,k\sep 1,\ldots,k)$ is a weight vector of weight
$-(\character_1+\cdots+\character_k)$;  the line through it is $B^-$ stable;
\item\mylabel{i.three}
$\det(n,\ldots,n-k+1\sep 1,\ldots,k)$ is a weight vector of weight
$-(\character_n+\cdots+\character_{n-k+1})$; the line through it is $B$-stable.
\end{enumerate}
\end{proposition}
\noindent
The proof of the proposition consists of elementary verifications.

We continue with the proof of the theorem.
It follows from (\ref{i.one}) of the proposition that 
\begin{eqnarray*}
\sectionp &:=&\det(1\sep1)\ \det(1,2\sep1,2)\ \cdots\ \det(1,\ldots,n-1\sep1,\ldots,n-1)\\
\sectionq &:=&\det(n\sep1)\ \det(n,n-1\sep1,2)\ \cdots\ \det(n,n-1,\ldots,2\sep1,\ldots,n-1)
\end{eqnarray*}
are both sections of $\lbundle(\rho)$.    We claim that the $(p-1)^\textup{st}$
power of $\sectionp\sectionq$ is the section of~$K^{1-p}$ that gives the $B$-canonical splitting
of~$G/B$.   By Theorem~\ref{t:canonical}, it is enough to check that this
section is a~$B$-canonical element, and for this we make use of the
characterization in~Proposition~\ref{p:canonical}.

Apply the proposition with $X=SL_n/B$.
As observed in~\S\ref{s.splittings} above, $\End_F(X)$
is naturally isomorphic to $\hnot(X,K^{1-p})$ (when $X$ is a non-singular
projective variety).  
By~Proposition~\ref{p.determinants}~(\ref{i.two}),
we can take $f_-$ to be $\sectionp^{p-1}$ 
(observe that
$\sum_{k=1}^{n-1}-(\character_1+\cdots+\character_k)=\sum_{k=1}^{n-1}-\fund_k=-\rho$,
and that the Steinberg module has lowest weight $-(p-1)\rho$);
by Proposition~\ref{p.determinants}~(\ref{i.three}), we can take $f_+$ 
to be $\sectionq^{p-1}$
(observe that
$\sum_{k=1}^{n-1}-(\character_n+\cdots+\character_{n-k+1})=\sum_{k=1}^{n-1}\weylinvolution(\character_1+\cdots+\character_k)$ where $\weylinvolution$ is the Weyl 
involution (whose action in the case of $SL_n$ on characters of~$T$ is given by
$\weylinvolution(\character_k)=-\character_{n-k+1})$, and so $\sectionq$ has
weight $\sum_{k=1}^{n-1}\weylinvolution(\fund_k)=\weylinvolution(\rho)=\rho$).
Observe also that
when global sections of line bundles on $G/B$ are identified
as functions on the group~$G$ that transform according to~(\ref{eq:transform1}),
the natural map $H^0(G/B,\lbundle_1)\tensor H^0(G/B,\lbundle_2)\to 
      H^0(G/B,\lbundle_1\tensor\lbundle_2)$ is just the ordinary multiplication
of functions.   This completes the proof that 
$(\sectionp\sectionq)^{p-1}$ is $B$-canonical.

Now fix a standard maximal parabolic $P$ of $SL_n$.   Given the claim
of the previous paragraph,  it is enough, in order to prove the theorem,
to show that the order of vanishing of 
the  section~$\sectionp\sectionq$ along~$P/B$ 
equals the codimension of $P/B$ in $SL_n/B$ (which is the same
as the dimension of~$SL_n/P$).     Since the identity coset
$eB$ (here $e$ denotes the identity element of the group~$SL_n$) belongs 
to~$P/B$ and $\sectionp$ does not vanish at $eB$,   
we have only to be concerned with~$\sectionq$.  By 
Proposition~\ref{p.determinants}~(\ref{i.three}), 
$\sectionq$ is a weight vector of weight $\rho$ 
for the action of $B$,
that is, $\sectionq(bg)=\rho(b)\cdot \sectionq(g)$ 
for $b\in B$ and $g\in SL_n$.
Thus the order of vanishing of~$\sectionq$ along $P/B$ 
is the same as that of~$\sectionq$
at the $T$-fixed point that is the `center of $P/B$', namely,
$w_0^PB$ (where $w_0^P$ denotes the longest element in the Weyl group of~$P$).

Let now $P$ be the maximal parabolic subgroup corresponding to $\fund_r$.
Then $w_0^P=(r,\ldots,1,n,\ldots,n-r+1)$.   The affine patch of~$SL_n/B$ centered
around~$w_0^PB$ given by $w_0^PB^-B$ consists of matrices of the following
explicit form:
\\[.1cm]
\begin{center}
\setlength{\unitlength}{.35cm}
\begin{picture}(15,15)(-1,-1)
\linethickness{.25mm}
\put(0,0){\line(0,1){13}}
\put(0,13){\line(1,0){.5}}
\put(0,0){\line(1,0){.5}}
\put(13,13){\line(1,0){.5}}
\put(13,0){\line(1,0){.5}}
\put(-1.1,0.25){$n$}
\put(-2.6,6.95){$r+1$}
\put(-1.1,8.25){$r$}\put(.5,8.25){1}
\multiput(1.5,9.25)(1,1){3}{$\cdot$}
\put(-1.1,12.25){$1$}\put(4.5,12.25){1}
\put(3.5,9.25){{\huge{\ensuremath{\star}}}}
\put(1,11.25){{\huge{0}}}
\put(2.5,3.5){{\huge{\ensuremath{\star}}}}
\put(5.75,0.25){$1$}\put(12.75,6.95){$1$}
\multiput(6.75,1.20)(1,0.96){6}{$\cdot$}
\put(8,5.25){{\huge{0}}}
\put(9.5,10){{\huge{0}}}
\put(10.5,2.25){{\huge{\ensuremath{\star}}}}
\put(13.5,0){\line(0,1){13}}
\put(.5,13.5){$1$}
\put(4.5,13.5){$r$}
\put(5.4,13.5){$r+1$}
\put(12.5,13.5){$n$}
\linethickness{.125mm}
\put(0.25,7.75){\line(1,0){13}}
\put(5.25,0){\line(0,1){13}}
\end{picture}
\end{center}
It is now elementary to check that 
$\det(n,\ldots,n-k+1\sep1,\ldots,k)$ 
vanishes at $w_0^P$ to order
\[
\textup{ord}_{P/B}\left(\det(n,\ldots,n-k+1\sep1,\ldots,k)\right)=
\left\{\begin{array}{ll}
      k & \textup{if $k\leq r$ and $k+r\leq n$} \\
      r & \textup{if $k\geq r$ and $k+r\leq n$} \\
      n-k & \textup{if $k\geq r$ and $k+r\geq n$} \\
      n-r & \textup{if $k\leq r$ and $k+r\geq n$} \\
          \end{array}\right.   \]
An elementary calculation using this shows
$\textup{ord}_{P/B}(\sectionq)=r(n-r)=\dim SL_n/P$.
This finishes the proof of the theorem in case~(1).

\mysubsection{The case of $G=Sp_{2m}$ and $P=P_n$}%
\mylabel{ss.sp2n}
Let $V$ be a $\field$-vector space with a non-degenerate skew-symmetric form
$\langle\ ,\ \rangle$.    The dimension of $V$ is then even, say~$2n$.
For $1\leq k\leq 2n$,  let $k^\ensuremath{\star}:=2n+1-k$.    Fix a basis 
$\basis_1,\ldots,\basis_{2n}$ of~$V$ such that 
\[
\langle \basis_i,\basis_j \rangle = \left\{
      \begin{array}{cl}
        1 & \textup{if $j=i^\ensuremath{\star}$ and $i<j$}\\
        -1 & \textup{if $j=i^\ensuremath{\star}$ and $i>j$}\\
        0 & \textup{if $j\neq i^\ensuremath{\star}$}
        \end{array}\right.
\]
and think of elements of $Sp_{2n}$ as $2n\times 2n$ matrices with 
respect to this basis (that preserve the form $\langle\ ,\ \rangle$).
The advantage of such a choice of basis is
this:  matrices in $Sp_{2n}$ that are diagonal form a maximal
torus (in~$Sp_{2n}$) and matrices that are upper triangular 
form a Borel subgroup.  We take $T$ and $B$ to be these.  We continue
to denote by $\character_k$ the restriction to $T$ of 
the character $\character_k$ of the diagonal torus of~$SL_{2n}$.
An easy verification shows that $\character_k=-\character_{k^\ensuremath{\star}}$.
The simple roots in order are $\character_1-\character_2$,
$\character_2-\character_3$, \ldots,
$\character_{n-1}-\character_{n}$, and
$\character_{n}-\character_{n+1}=
2\character_n$,  and the corresponding fundamental weights are
\[
\fund_1=\character_1,\quad
\fund_2=\character_1+\character_2,\quad
\ldots,\quad
\fund_n=\character_1+\cdots+\character_n.\]

Consider the following functions on $Sp_{2n}$ (the symbols on the
right hand side denote functions on $SL_{2n}$ in the notation
defined in~\S\ref{ss.sln} above, 
and now they also denote the restriction to $Sp_{2n}$
of those functions):
\begin{eqnarray*}
\sectionp&:=& \det(1\sep1)\,\det(1,2\sep1,2)\,\cdots\,
       \det(1,2,\ldots,n\sep 1,2,\ldots,n)\\
\sectionq & := & \det(2n\sep1)\,\det(2n,2n-1\sep1,2)\,\cdots\,
       \det(2n,2n-1,\ldots,n\sep 1,2,\ldots,n)    
\end{eqnarray*}
Just as in the case of $SL_n$,  an easy verification using 
Proposition~\ref{p.determinants} shows that $\sectionp$, $\sectionq$
are sections of $\lbundle(\rho)$  and that the $(p-1)^\textup{st}$
power of $\sectionp\sectionq$ is the section of $K^{1-p}$ that
gives the $B$-canonical splitting of $Sp_{2n}/B$.  

Let $P=P_n$ be the maximal parabolic subgroup corresponding to
the fundamental weight $\fund_n=\character_1+\cdots+\character_n$.
We calculate the order of vanishing of $\sectionp\sectionq$ along
$P/B$ and show that it equals the codimension of $P/B$ in $Sp_{2n}/B$ 
(which equals $\sum_{k=1}^n k= {{n+1}\choose{2}}$).
This calculation too runs parallel to the case of $SL_n$.    Just
as in that case, we reduce to considering the order of vanishing
of $\sectionq$ at the point $w_0^PB$ of $P/B$.

The affine patch around $w_0^PB$ of $Sp_{2n}/B$ given by $w_0^PB^-B$ consists
of matrices having the following explicit form:
\\[2mm]\begin{center}
\hspace{2cm}
\begin{minipage}{4.2cm}
\setlength{\unitlength}{.35cm}
\begin{picture}(12,12)(2,2)
\linethickness{.25mm}
\put(0,3){\line(0,1){10}}
\put(10.5,3){\line(0,1){10}}
\put(0,13){\line(1,0){.5}}
\put(0,3){\line(1,0){.5}}
\put(10,13){\line(1,0){.5}}
\put(10,3){\line(1,0){.5}}
\put(-1.5,3.25){$2n$}
\put(-2.6,6.95){$n+1$}
\put(-1.1,8.25){$n$}\put(.5,8.25){1}
\multiput(1.5,9.25)(1,1){3}{$\cdot$}
\put(-1.1,12.25){$1$}\put(4.5,12.25){1}
\put(3.5,9.25){{\huge{\ensuremath{\star}}}}
\put(1,11.25){{\huge{0}}}
\put(2.5,5){{\huge{\ensuremath{\star}}}}
\put(5.75,3.25){$1$}\put(9.75,6.95){$1$}
\multiput(6.75,4.20)(1,0.96){3}{$\cdot$}
\put(6.5,6.25){{\huge{0}}}
\put(7.5,10){{\huge{0}}}
\put(8.5,4.25){{\huge{\ensuremath{\star}}}}
\put(.5,13.5){$1$}
\put(4.5,13.5){$n$}
\put(5.4,13.5){$n+1$}
\put(9.25,13.5){$2n$}
\linethickness{.125mm}
\put(0.25,7.75){\line(1,0){10}}
\put(5.25,3){\line(0,1){10}}
\end{picture}\end{minipage}
\begin{minipage}{3.5cm}
{\raggedright
The inderminates in positions~$\star$ are
not algebraically independent.}
\end{minipage}
\end{center}

We claim that the order of vanishing of the section
$\det(2n,\ldots,2n-k+1\sep1,\ldots,k)$
(for $1\leq k\leq n$) at $w_0^PB$ is $k$.
The order of vanishing being $k$ for $SL_{2n}$,  it follows
immediately that now it is no less than~$k$.   However,  since,
unlike in the case of $SL_{2n}$,  the entries indicated by~\ensuremath{\star} 
in the matrix above are not algebraically independent,  it requires
some proof that it is no more than~$k$.    For this,
we specialize:  set all the variables~\ensuremath{\star} equal to~$0$ except those
on the anti-diagonal in the~$n\times n$ matrix in the bottom left
corner and take those on the anti-diagonal to be algebraically
independent variables (this is a valid specialization in the sense that the
resulting matrices are inside $Sp_{2n}$).    The restriction of
$\det(2n,\ldots,2n-k+1\sep 1,\ldots,k)$ to this closed set (which
is an affine $n$-space) is given by the product of $k$ indeterminates
and so has order of vanishing exactly~$k$ at the origin.    This finishes
the proof of the theorem in case~(2).

\mysubsection{The case $G=SO_{2n}$ and $P=P_n$}%
\mylabel{ss.so2n}
This proof runs mostly parallel to that for~$Sp_{2n}/P_n$,
but there are two notable differences---firstly, the number of
factors in the definitions below of the
functions~$\sectionp$ and~$\sectionq$ does not equal the number of
fundamental weights;  secondly,  the specialization when $n$ is odd
does not work in the same fashion as when it is even.

Let $V$ be a $\field$-vector space of even dimension,  say~$2n$, 
with a non-degenerate symmetric form
$\left(\ ,\ \right)$.    
For $1\leq k\leq 2n$,  let $k^\ensuremath{\star}:=2n+1-k$.    Fix a basis 
$\basis_1,\ldots,\basis_{2n}$ of~$V$ such that 
\[
\left( \basis_i,\basis_j \right) = \left\{
      \begin{array}{cl}
        1 & \textup{if $j=i^\ensuremath{\star}$}\\
                0 & \textup{if $j\neq i^\ensuremath{\star}$}
        \end{array}\right.
\]
and think of elements of $SO_{2n}$ as $2n\times 2n$ matrices with 
respect to this basis (that preserve the form $\left(\ ,\ \right)$
and have determinant~$1$).
The advantage of such a choice of basis is
this:  matrices in $SO_{2n}$ that are diagonal form a maximal
torus (in~$SO_{2n}$) and matrices that are upper triangular 
form a Borel subgroup.  We take $T$ and $B$ to be these.  We continue
to denote by $\character_k$ the restriction to $T$ of 
the character $\character_k$ of the diagonal torus of~$SL_{2n}$.
An easy verification shows that $\character_k=-\character_{k^\ensuremath{\star}}$.
The simple roots in order are $\character_1-\character_2$,
$\character_2-\character_3$, \ldots,
$\character_{n-1}-\character_{n}$, and
$\character_{n-1}-\character_{n+1}=\character_{n-1}+\character_n$,
and the corresponding fundamental weights are
\begin{multline*}
\fund_1=\character_1,\quad
\fund_2=\character_1+\character_2,\quad
\ldots,\quad
\fund_{n-2}=\character_1+\cdots+\character_{n-2},\\[2mm]
\fund_{n-1}=\frac{1}{2}\left(\character_1+\cdots+\character_{n-1}%
      -\character_n\right),\quad
\fund_n=\frac{1}{2}\left(\character_1+\cdots+\character_{n-1}+
                                               \character_n\right).
\end{multline*}

Consider the following functions on $SO_{2n}$ (the symbols on the
right hand side denote functions on $SL_{2n}$ in the notation defined
in~\S\ref{ss.sln}, and now they represent the restriction to $SO_{2n}$ of those
functions):
\begin{eqnarray*}
\sectionp \! &:=&\! \det(1\sep1)\,\det(1,2\sep1,2)\,\cdots\,
       \det(1,2,\ldots,n-1\sep 1,2,\ldots,n-1)\\
\sectionq\! & := & \!\det(2n\sep1)\,\det(2n,2n-1\sep1,2)\,\cdots\,
       \det(2n,2n-1,\ldots,n+1\sep 1,2,\ldots,n-1)    
\end{eqnarray*}
Note that the weight $\character_1+\cdots+\character_{n-1}$ of the
$(n-1)^\textup{st}$ factor of $\sectionq$ 
for the action of~$T$ is precisely the
sum of $\fund_{n-1}$ and $\fund_n$  (the weight of $\sectionp$ is
the negative of that of~$\sectionq$).   Using this and 
Proposition~\ref{p.determinants},  it follows,
just as in the case of~$SL_n$, that $\sectionp$, $\sectionq$
are sections of $\lbundle(\rho)$  and that the $(p-1)^\textup{st}$
power of $\sectionp\sectionq$ is the section of $K^{1-p}$ that
gives the $B$-canonical splitting of~$SO_{2n}/B$.  

Let $P=P_n$ be the maximal parabolic subgroup corresponding to the
fundamental weight $\fund_n=\frac{1}{2}(\character_1+\cdots+\character_n)$.  We
calculate the order of vanishing of $\sectionp\sectionq$ along $P/B$
and show that it equals the codimension of $P/B$ in $SO_{2n}/B$ (which
is $\sum_{k=1}^{n-1}k={n\choose 2}$).  This
calculation too runs parallel to the case of $SL_n$.  Just as in that
case, we reduce to considering the order of vanishing of $\sectionq$
at the point $w_0^PB$ of $P/B$.

The affine patch around $w_0^PB$ of $SO_{2n}/B$ given by $w_{0}^{P}B^-B$ consists
of matrices having the following explicit form:
\\[2mm]\begin{center}
\hspace{2cm}
\begin{minipage}{4.2cm}
\setlength{\unitlength}{.35cm}
\begin{picture}(12,12)(2,2)
\linethickness{.25mm}
\put(0,3){\line(0,1){10}}
\put(10.5,3){\line(0,1){10}}
\put(0,13){\line(1,0){.5}}
\put(0,3){\line(1,0){.5}}
\put(10,13){\line(1,0){.5}}
\put(10,3){\line(1,0){.5}}
\put(-1.5,3.25){$2n$}
\put(-2.6,6.95){$n+1$}
\put(-1.1,8.25){$n$}\put(.5,8.25){1}
\multiput(1.5,9.25)(1,1){3}{$\cdot$}
\put(-1.1,12.25){$1$}\put(4.5,12.25){1}
\put(3.5,9.25){{\huge{\ensuremath{\star}}}}
\put(1,11.25){{\huge{0}}}
\put(2.5,5){{\huge{\ensuremath{\star}}}}
\put(5.75,3.25){$1$}\put(9.75,6.95){$1$}
\multiput(6.75,4.20)(1,0.96){3}{$\cdot$}
\put(6.5,6.25){{\huge{0}}}
\put(7.5,10){{\huge{0}}}
\put(8.5,4.25){{\huge{\ensuremath{\star}}}}
\put(.5,13.5){$1$}
\put(4.5,13.5){$n$}
\put(5.4,13.5){$n+1$}
\put(9.25,13.5){$2n$}
\linethickness{.125mm}
\put(0.25,7.75){\line(1,0){10}}
\put(5.25,3){\line(0,1){10}}
\end{picture}\end{minipage}
\begin{minipage}{3.5cm}
{\raggedright
The inderminates in positions~$\star$ are
not algebraically independent.}
\end{minipage}
\end{center}

Unlike in the case of $SL_n$,  the entries indicated by~\ensuremath{\star} in the
picture above are not algebraically independent.  
Since $w_0^PB$ belongs to $SO_{2n}/B$,  it follows that the 
ideal defining the affine patch of $SO_{2n}/B$ as a closed subset 
of $SL_{2n}/B$ is contained in the maximal ideal generated by the
indeterminates~\ensuremath{\star}.     We can therefore conclude that the order
of vanishing of $\det(2n,\ldots,2n-k+1\sep1,\ldots,k)$ at $w_0^PB$
is no less than~$k$,  so that the order of vanishing of~$\sectionq$
at $w_0^PB$ is no less than $\sum_{k=1}^{n-1} k=(n-1)n/2=\dim SO_{2n}/P$.
To see that the order of vanishing equals this lower bound,
we specialize.   We first do this when $n$ is even.
Set all the variables~\ensuremath{\star} equal to~$0$ except those
on the anti-diagonal in the~$n\times n$ matrix in the bottom left
corner;   set the variables~\ensuremath{\star} that are on the anti-diagonal on
rows $3n/2$ through $2n$ to be algebraically independent variables,
say $X_{3n/2}$,~\ldots,~$X_{2n}$; and set the variables~\ensuremath{\star} that are on 
the anti-diagonal on rows $n+1$ through~$3n/2-1$ to be 
$-X_{2n}$,~\ldots,~$-X_{3n/2}$ (this is a valid specialization 
in the sense that the
resulting matrices are inside $SO_{2n}$).    The restriction of
$\det(2n,\ldots,2n-k+1\sep 1,\ldots,k)$ to this closed set (which
is an affine $n/2$-space) is given by the product of $k$ (possibly
repeated and with signs) indeterminates
and so has order of vanishing exactly~$k$ at the origin.    This finishes
the proof of the theorem in case~(3) when $n$ is even.

Now suppose that $n$ is odd.   Specialize as follows:  set equal
to~$0$ all variables~$\star$ 
not in the $n\times n$ matrix in the lower left hand
corner;  take the~$n\times n$ matrix in the lower left hand corner to be
a generic skew-symmetric matrix.     This is a valid specialization.
The resulting space is an affine $(n-1)n/2$-space.
Our goal is to show that 
the restriction of $\det(2n,\ldots,2n-k+1\sep 1,\ldots,k)$ 
to this affine space does not lie in the $(k+1)^\textup{st}$ power
of the maximal ideal corresponding to the origin, and this will
be reached once we prove the following claim:
let $V$ be an $n$-dimensional vector space
with a skew-symmetric form~$\langle\ ,\ \rangle$ 
of rank $n-1$ (such a form exists);
for $k$ an integer, $1\leq k\leq n-1$, there exist vectors $e_1,\ldots,e_k$
and $e_{n-k+1},\ldots,e_n$ such that 
the matrix $\left(\langle e_i,e_j\rangle\right)$,
$n-k+1\leq j\leq n$ and $1\leq i\leq k$,  is invertible.
(To see why it suffices to prove the claim,  we think of the skew-symmetric
$n\times n$ matrix in the bottom left hand corner as defining the form
with respect to some basis.)    

To prove the claim,  let~$W$ be a $k$-dimensional subspace of~$V$ that does 
not meet the radical\footnote{Recall that the radical consists
of the elements $w$ of $W$ such that 
$\langle w,x\rangle=0$ $\forall$ $x\in W$.} 
(which is $1$ dimensional by our hypothesis).   
Then $W^\perp$ has dimension~$n-k$.  Let~$W'$
to be a $k$-dimensional subspace that meets $W^\perp$ trivially and
intersects~$W$ in only a subspace of dimension~$2k-n$.   Choose
$e_{n-k+1},\ldots,e_k$ to be a basis of~$W\cap W'$, and extend
it to a basis $e_1,\ldots,e_k$ for $W$ and to a basis $e_{n-k+1},\ldots,e_n$
of $W'$. \end{myproof}
\mysection{Conclusion: Wahl's conjecture holds for the $G/P$ of 
Theorem~\ref{t.main}}\mylabel{s.conclusion}
We assume in this section that the characteristic~$p$ is odd.
In order to prove Wahl's conjecture for~$G/P$,  it is enough,
by Theorem~\ref{t.lmp},  that there exist a splitting of $G/P\times G/P$ which
compatibly splits the diagonal copy of~$G/P$ with maximal multiplicity.
We now argue that such a splitting exists for the $G/P$ as in
Theorem~\ref{t.main}.     In fact, we show that the canonical splitting
of $G/P\times G/P$ defined in~\S\ref{s.canonical} above has the desired
property.   Our argument follows that in~\cite{mp}.

First recall the following explicit description of the canonical
splitting of $G/B\times G/B$.
Let~$D$ be the union of the Schubert divisors in~$G/B$ and~$\tilde{D}$
the union of the opposite Schubert divisors.   Let $p_1$ and~$p_2$ denote
respectively the first and second projections of $G/B\times G/B$ onto
$G/B$.     Then, as shown in~\cite{mr2}, $\sigma^{p-1}$ 
is a splitting of $G/B\times G/B$ (that compatibly splits the diagonal
copy of~$G/B$),  where $\sigma$ denotes the section of the 
the canonical bundle of~$G/B\times G/B$ given by 
the divisor~$p_1^*D+G\times^BD+p_2^*\tilde{D}$.  The presence of~$p_1^*D$
means that this splitting splits compatibly $X\times G/B$ for $X$ any
Schubert variety in~$G/B$;   in particular,  $e\times G/B$ is compatibly
split.     Thus, in order show that~$\sigma^{p-1}$ is the canonical
splitting, 
it suffices, by Theorem~\ref{t.canonical.2},   
to show that $\sigma^{p-1}$ is $B$-canonical.

Before turning to the proof of this,  
let us finish the proof of Wahl's conjecture.
It has been proved in~\S\ref{s.theorem} that the order of vanishing
of~$D$ along~$P/B$ equals the codimension of~$P/B$ (in $G/B$).
The occurrence of $G\times^B D$ in the divisor defining~$\sigma$
and the smoothness of the natural map $G/B\times G/B\to G/P\times G/P$ 
along $G\times^B P/B$ together imply (by~\cite[Lemma~3.1]{mp}
or~\cite[1.3.E.13]{bk}) that the image of $G\times^B P/B$ is compatibly
split with maximal multiplicity (under the induced splitting of 
$G/P\times G/P$).   But the image of $G\times^B P/B$ is the diagonal.
This completes the proof of Wahl's conjecture.

To prove that~$\sigma^{p-1}$ is $B$-canonical,
we use the characterization in Proposition~\ref{p:canonical}.  
As observed in~\S\ref{s.splittings},   we have
\[ \End_F(G/B\times G/B)\,\,\cong\,\, H^0(G/B\times G/B, 
                    \lbundle(2(p-1)\rho,2(p-1)\rho))    \]
The right hand side above is naturally isomorphic to
\[
H^0(G_1/B_1,\lbundle(2(p-1)\rho_1))
\,\,\tensor\,\, H^0(G_2/B_2,\lbundle(2(p-1)\rho_2))\]
where the subscripts~$1$
and~$2$ are used to denote respectively
objects associated to the copies~$G\times e$
and $e\times G$ of~$G$ in~$G\times G$.   And there are natural maps
\[ \steinberg_1\tensor\steinberg_1\,\to\, H^0(G_1/B_1,\lbundle(2(p-1)\rho_1)
\quad\quad\steinberg_2\tensor\steinberg_2\,\to\, H^0(G_2/B_2,\lbundle(2(p-1)\rho_2)\]
So we have a natural map
\[ 
     \steinberg_1\tensor    \steinberg_1\tensor   
            \steinberg_2\tensor    \steinberg_2\,\, \longrightarrow\,\,
\End_F(G/B\times G/B) \]

As~$G$-modules,  we have $\steinberg_1\cong\steinberg\cong\steinberg_2$.
And,  since~$\steinberg$ is irreducible and self-dual,   there exists a
unique $G$-invariant element---call it $\sectionr$---in $\steinberg_1\tensor
\steinberg_2$.   On the other hand,  as seen in~\S\ref{s.theorem},  the
sections~$\sectionp^{p-1}$ and~$\sectionq^{p-1}$ (with $\sectionp$ and
$\sectionq$ being defined as there) can respectively 
be taken to be a (non-zero) lowest weight vector in~$\steinberg_2$
and a (non-zero) highest weight vector in~$\steinberg_1$.

Putting all the above together,  we get a $\field$-linear $B$-module map
\[
        \field_{(p-1)\rho}\tensor \steinberg \longrightarrow
\steinberg_1\tensor\steinberg_1\tensor\steinberg_2\tensor\steinberg_2
                                                                \]
sending~$f_+\tensor f_-$ to $\sectionq^{p-1}\tensor\sectionr\tensor
\sectionp^{p-1}$.   But the splitting section~$\sigma^{p-1}$ is just the
natural image in~$\End_F(G/B\times G/B)$ of 
$\sectionq^{p-1}\tensor\sectionr\tensor\sectionp^{p-1}$:
$p_1^*D$ corresponds to $\sectionq^{p-1}$,
$\sectionr$ to $G\times^B D$,
and $p_2^*\tilde{D}$ to $\sectionp^{p-1}$.  This finishes the
proof that~$\sigma^{p-1}$ is $B$-canonical.

	\newcommand\citenumfont[1]{\textbf{#1}}
\ifthenelse{\equal{\finalized}{no}}{
\bibliographystyle{bibsty-final-no-issn-isbn}
\bibliography{abbrev,ref}

\begin{thebibliography}{10}
\expandafter\ifx\csname url\endcsname\relax
  \def\url#1{{\tt #1}}\fi
\expandafter\ifx\csname urlprefix\endcsname\relax\def\urlprefix{URL }\fi

\bibitem{bour}
N.~Bourbaki, {\em \'{E}l\'ements de math\'ematique. {F}asc. {XXXIV}. {G}roupes
  et alg\`ebres de {L}ie. {C}hapitre {IV}: {G}roupes de {C}oxeter et syst\`emes
  de {T}its. {C}hapitre {V}: {G}roupes engendr\'es par des r\'eflexions.
  {C}hapitre {VI}: syst\`emes de racines\/}, Actualit\'es Scientifiques et
  Industrielles, No. 1337, Hermann, Paris, 1968.

\bibitem{bk}
M.~Brion and S.~Kumar, {\em Frobenius splitting methods in geometry and
  representation theory\/}, vol. 231 of Progress in Mathematics, Birkh\"auser
  Boston Inc., Boston, MA, 2005.

\bibitem{kumar}
S.~Kumar, {\em Proof of {W}ahl's conjecture on surjectivity of the {G}aussian
  map for flag varieties\/}, Amer. J. Math., {\bf 114}, no.~6, 1992,
  pp.~1201--1220.

\bibitem{lmp}
V.~Lakshmibai, V.~B. Mehta, and A.~J. Parameswaran, {\em Frobenius splittings
  and blow-ups\/}, J. Algebra, {\bf 208}, no.~1, 1998, pp.~101--128.

\bibitem{mat}
O.~Mathieu, {\em Filtrations of {$G$}-modules\/}, Ann. Sci. \'Ecole Norm. Sup.
  (4), {\bf 23}, no.~4, 1990, pp.~625--644.

\bibitem{mat00}
O.~Mathieu, {\em Tilting modules and their applications\/}, in: {\em Analysis
  on homogeneous spaces and representation theory of Lie groups, Okayama--Kyoto
  (1997)\/}, vol.~26 of Adv. Stud. Pure Math., Math. Soc. Japan, Tokyo, 2000,
  pp. 145--212.

\bibitem{mp}
V.~B. Mehta and A.~J. Parameswaran, {\em On {W}ahl's conjecture for the
  {G}rassmannians in positive characteristic\/}, Internat. J. Math., {\bf 8},
  no.~4, 1997, pp.~495--498.

\bibitem{mr}
V.~B. Mehta and A.~Ramanathan, {\em Frobenius splitting and cohomology
  vanishing for {S}chubert varieties\/}, Ann. of Math. (2), {\bf 122}, no.~1,
  1985, pp.~27--40.

\bibitem{mr2}
V.~B. Mehta and A.~Ramanathan, {\em Schubert varieties in {$G/B\times G/B$}\/},
  Compositio Math., {\bf 67}, no.~3, 1988, pp.~355--358.

\bibitem{vdk}
W.~van~der Kallen, {\em Lectures on {F}robenius splittings and
  {$B$}-modules\/}, Published for the Tata Institute of Fundamental Research,
  Bombay, by Springer-Verlag, Berlin, 1993. Notes by S. P. Inamdar.

\bibitem{wahl}
J.~Wahl, {\em Gaussian maps and tensor products of irreducible
  representations\/}, Manuscripta Math., {\bf 73}, no.~3, 1991, pp.~229--259.

\end{thebibliography}
}{
}

\end{document}